\documentclass{article}
\usepackage{amsfonts}
\usepackage{amsmath}
\usepackage{footnote}
\usepackage{multicol}
\usepackage{booktabs}
\usepackage[flushleft]{threeparttable}
\usepackage[font=small,labelfont=bf]{caption}

\usepackage{graphicx}
\usepackage{amsthm}

\newtheorem{definition}{Definition}
\newtheorem{remark}{Remark}
\newtheorem{proposition}{Proposition}
\newtheorem{theorem}{Theorem}
\newtheorem{corollary}{Corollary}
\newtheorem{example}{Example}

\begin{document}
%


\title{Continuity up to a covering and connectedness}

\author{Emin Durmishi\\ {\em University of Tetova, North Macedonia} \\ {\em email: emin.durmishi@unite.edu.mk}}

\maketitle
\begin{abstract}
One of the ways that connectedness has been studied through the history of topology is by using chains, the so called chain connectedness. Here we combine this notion together with continuity up to a covering to provide the inheritance of connectedness for the topological spaces even when there is no continuous surjection between them.
\end{abstract}
\textbf{Keaywords:} Continuity up to a covering, chain connectedness, connectedness.        



\section{Introduction}

On the second half of the twentieth century there have been introduced many generalizations of the concept of continuity between two topological spaces. For some of the definitions, their properties and their comparison we refer to \cite{[6],[7],[8]} and their references. Here we will focus on the notion of continuity up to a covering.

The notion of $\varepsilon$-continuity (see \cite{[5]}) led to a more general definition, by Kieboom (see \cite{[11]}), that of continuity up to some covering $\mathcal{V}$, i.e., $\mathcal{V}$-continuity. Proximate nets, used in intrinsic shape theory, are defined by using the notion of $\mathcal{V}$-continuity (see \cite{[11],[1]}). 

Connectedness has first been attempted to be defined, by Cantor, by using the notion of chain. In \cite{[2]} the notion of chain connectedness in a topological space is presented, and in \cite{[3]} and \cite{[4]} this notion is generalized to subsets of the space relative to the space.

Both notions use open coverings on their definitions.

\subsection{$\mathcal{V}$-continuous functions}

Through the text, by covering it will always be meant an open covering.

Let $f:X\to Y$ be a function between two topological spaces and let $\mathcal{V}$ be a covering for $Y$.

\begin{definition}
The function $f$ is said to be {\em $\mathcal{V}$-continuous} if for any $x\in X$ there is a neighbourhood $U$ of $x$ and $V\in \mathcal{V}$ such that $f(U)\subseteq V$.
\end{definition}

This is equivalent of finding a covering $\mathcal{U}$ of $X$ such that for every $x\in X$ there is a neighbourhood $U\in \mathcal{U}$ of $x$ and a $V\in \mathcal{V}$ for which $f(U)\subseteq V$.

\begin{remark}\label{rem}
Notice that if the covering $\mathcal{V}$ contains $Y$, then any function $f:X\to Y$ is $\mathcal{V}$-continuous.
\end{remark}

For two families of sets $\mathcal{U}$ and $\mathcal{V}$ we say that {\em $\mathcal{U}$ refines $\mathcal{V}$} if for any $U\in \mathcal{U}$ there exists a $V\in \mathcal{V}$ such that $U\subseteq V$, and we denote it by $\mathcal{U} \prec \mathcal{V}$.

If $f:X\to Y$ is a $\mathcal{V}$-continuous function and $\mathcal{W}$ is a covering of $Y$ such that $\mathcal{V} \prec \mathcal{W}$, then $f$ is $\mathcal{W}$-continuous.

\begin{proposition}
If for every $x\in X$ there exists a $V\in \mathcal{V}$ and a subset $W\subseteq V$ (not necessarily open) such that $f(x)\in W$ and $f^{-1}(W)$ is open in $X$, then $f:X\to Y$ is a $\mathcal{V}$-continuous function.
\end{proposition}

\begin{proof}
Let $f^{-1}(W)=U$. Since $f(x)\in W$, then $x\in f^{-1}(W)$, hence $U$ is a neighbourhood of $x$ such that
$$f(U)=f(f^{-1}(W))\subseteq W\subseteq V.$$
\end{proof}

\begin{proposition}
If $f:X\to Y$ is $\mathcal{V}$-continuous, $g:Y\to Z$ is $\mathcal{W}$-continuous and $g\left( \mathcal{V} \right) \prec \mathcal{W}$, then $g\circ f:X\to Z$ is $\mathcal{W}$-continuous.
\end{proposition}

It is obvious that continuity implies continuity up to a covering for any possible covering. The following theorem is a generalization of Proposition 1.3 (ii) from \cite{[11]}.

\begin{theorem}
Let $Y$ be a $T_1$ space. The function $f:X\to Y$ is continuous if and only if it is $\mathcal{V}$-continuous for every covering $\mathcal{V}$ of $Y$.
\end{theorem}

\begin{proof}
The necessary condition is trivial.

Suppose that $f:X\to Y$ is $\mathcal{V}$-continuous for any covering $\mathcal{V}$ of $Y$. Let $V\subseteq Y$ be an arbitrary open set and let $x\in f^{-1}(V)$. Then $f(x)\in V$ and, since $Y$ is $T_1$, $\forall y\in Y, \; y\neq f(x)$, there exists a neighbourhood $V_y$ of $y$ such that $f(x)\notin V_y$. Let $\mathcal{V}=\{ V_y\vert y\in Y, \; y\neq f(x) \} \cup \{ V \}$. Then $\mathcal{V}$ is a covering for $Y$, and since the only element of $\mathcal{V}$ that contains $f(x)$ is $V$, and on the other hand, by the assumption, $f$ is $\mathcal{V}$-continuous, there must exists a neighbourhood $U$ of $x$ such that $f(U)\subseteq V$, i.e., $x\in U\subseteq f^{-1}(V)$, therefore $f^{-1}(V)$ is open in $X$.
\end{proof}

The following example shows that $T_1$ is a crucial condition for the sufficiency of the theorem.

\begin{example}
Let $X=\mathbb{R}$ with the standard topology and let $Y=\{ a,b\}$ equipped with the topology $\tau_Y=\{ \emptyset , \{ a\} ,Y \}$. Then the function $f:X\to Y$ defined by $$f(x)=\begin{cases} a, & \textrm{ if } x\geq 0, \\ b, & \textrm{ if } x< 0 \end{cases}$$
is $\mathcal{V}$-continuous for any covering $\mathcal{V}$ of $Y$, but it is not continuous at $x=0$.
\end{example}

If $X$ is a discrete space, then any function $f:X\to Y$ is continuous, hence $\mathcal{V}$-continuous for any covering $\mathcal{V}$ of $Y$.

\subsection{Chain connectedness}

Connectedness has been attempted to be described by using the notion of chain since the time of Cantor. In \cite{[2]} connectedness of the topological space is characterized by using the chains for each covering, and further studied in \cite{[4]} and \cite{[3]}. Statements presented in this subsection are reformulations from these articles. 

Let $x$ and $y$ be two points in the topological space $X$ and $\mathcal{U}$ be a covering of $X$.

\begin{definition}
A {\em chain in $\mathcal{U}$ that joins $x$ and $y$} is a finite sequence $U_1,U_2,\ldots ,U_n$ of elements in $\mathcal{U}$ such that $x\in U_1, \; y\in U_n$ and $U_i\cap U_{i+1}\neq \emptyset , \forall i\in \{ 1,2,\ldots ,n-1 \} .$

If such a chain exists, we say that {\em $x$ and $y$ are $\mathcal{U}$-chain connected} and denote that with $x\mathop\sim\limits_{\mathcal{U}} y.$

If for every covering of the space we can find a chain joining $x$ and $y$ we say that {\em $x$ and $y$ are chain connected in $X$} and denote that with $x\sim y$.
\end{definition}

Both $\mathop\sim\limits_{\mathcal{U}}$ and $\sim$ are equivalence relations on $X$. Classes of equivalences of $\mathop\sim\limits_{\mathcal{U}}$ consists of clopen sets, while those of $\sim$ coincide with quasicomponents of the space $X$.

We introduce the notion of $\mathcal{U}$-chain connected sets.

\begin{definition}
A subset $C$ of the topological space $X$ is said to be {\em $\mathcal{U}$-chain connected in $X$} if every two points of $C$ are $\mathcal{U}-$chain connected in $X$.
\end{definition} 

In \cite{[4]} and \cite{[3]} the notion of chain connected set relative to the space is introduced.

\begin{definition}
A subset $C$ of the topological space $X$ is said to be {\em chain connected in $X$} if every two points of $C$ are chain connected in $X$.
\end{definition}

For a point $x\in X$, the maximal $\mathcal{U}$-chain connected set containing $x$ is called the {\em $\mathcal{U}$-chain component of $x$} and it is practically the class of equivalence of the $\mathcal{U}$-chain relation which contains $x$, while the maximal chain connected set containing $x$ is called the {\em chain component of $x$} and it is practically the class of equivalence of the chain relation which contains $x$.

Obviously, a set is chain connected in $X$ if and only if it is $\mathcal{U}$-chain connected for every covering $\mathcal{U}$ of $X$.

Any connected set is also chain connected in the related topological space, but chain connectedness does not imply connectedness (for counterexamples we refer to \cite{[3],[4]}).

However, the following statement, characterizes the connectedness of topological spaces by using chain connectedness.

\begin{theorem}
The topological space $X$ is connected if and only if it is chain connected in $X$.
\end{theorem}

By the following theorem we show that chain connectedness is inherited by continuous functions.

\begin{theorem}
If $f:X\to Y$ is a continuous function and $C\subseteq X$ is chain connected in $X$, then $f(C)$ is chain connected in $Y$. 
\end{theorem}

Hence, chain connectedness is a topological invariance, as stated below.

\begin{corollary}
If $f:X\to Y$ is a homeomorphism, then $C\subseteq X$ is chain connected in $X$ if and only if $f(C)$ is chain connected in $Y$.
\end{corollary}

In \cite{[9]} it is proved the following statement considering the product of chain connected sets in the product space.

\begin{theorem}
If $C_i$ are chain connected sets in $X_i, \; \forall i\in I$, then $\prod\limits_{i\in I}C_i$ is a chain connected set in $\prod\limits_{i\in I}X_i$ equipped with the product topology.
\end{theorem}

\section{Continuity up to a covering, chain connectedness and connectedness}

Since both continuity up to a covering and chain connectedness are defined by using coverings of space, it is expected to have new results by combining these notions. In this section we show the inheritance of chain connectedness by using continuous up to a covering functions. Then, we use these results to proof some results for connectedness of topological spaces. We generalize the well known fact that the connectedness is inherited by a continuous surjection.

\subsection{Continuity up to a covering and chain connectedness}

Let $\mathcal{V}$ be a covering of the topological space $Y$, $f:X\to Y$ be a $\mathcal{V}$-continuous function and $\mathcal{U}$ be a covering of the topological space $X$ such that for every $U\in \mathcal{U}$ there is a $V\in \mathcal{V}$ for which $f(U)\subseteq V$.

\begin{proposition}\label{prop11}
If $X$ is a $\mathcal{U}$-chain connected space, then $f(X)$ is a $\mathcal{V}$-chain connected subspace of $Y$.
\end{proposition}

\begin{proof}
Let $x$ and $y$ be two arbitrary points in $X$ and $U_1,U_2,\ldots ,U_n$ be a chain in $\mathcal{U}$ joining $x$ and $y$. Then there exists a sequence $V_1,V_2,\ldots ,V_n$ in $\mathcal{V}$ such that $f(U_i)\subseteq V_i, \; \forall i\in \{ 1,2,\ldots ,n \}$. Therefore, for two arbitrary points $f(x),f(y)\in f(X)$ the sequence $V_1\cap f(X),V_2 \cap f(X),\ldots , V_n\cap f(X)$ is a chain in $\mathcal{V}_{f(X)}=\left\{ V\cap f(X) \vert V\in \mathcal{V} \right\}$ that joins them.
\end{proof}

\begin{corollary}
If $X$ is a connected topological space and $f:X\to Y$ is a $\mathcal{V}$-continuous function, then $f(X)$ is a $\mathcal{V}$-chain connected subspace of $Y$.
\end{corollary}

\begin{proof}
Since $X$ is connected, from the Proposition 1.1 in \cite{[2]}, for any $x,y\in X$ and any covering $\mathcal{U}$ of $X$ there is a chain in $\mathcal{U}$ that joins $x$ and $y$. Since $f$ is $\mathcal{V}$-continuous, we can chose a covering $\mathcal{U}$ of $X$ as in the conditions of Proposition \ref{prop11}, therefore, from the same proposition, $f(X)$ is a $\mathcal{V}$-chain connected subspace.
\end{proof}

\begin{corollary}
If $X$ is a connected topological space and there exists a $\mathcal{V}$-continuous surjection $f:X\to Y$, then $Y$ is $\mathcal{V}$-chain connected.
\end{corollary}

\subsection{Continuity up to a covering and connectedness}

The following claim proves the inheritance of connectedness by using the existence of a $\mathcal{V}$-continuous function for each covering $\mathcal{V}$.

\begin{theorem}\label{lasttheorem}
If $X$ is a connected topological space and for each covering $\mathcal{V}$ of $Y$ there exists a $\mathcal{V}$-continuous surjection $f_{\mathcal{V}}:X\to Y$, then $Y$ is connected.
\end{theorem}

\begin{proof}
Let $\mathcal{V}$ be an arbitrary covering of $Y$ and let $f_{\mathcal{V}}(x),f_{\mathcal{V}}(y)\in Y$. Since $f_{\mathcal{V}}$ is $\mathcal{V}$-continuous, there exists a covering $\mathcal{U}$ of $X$ such that for each $U\in \mathcal{U}$ there exists $V\in \mathcal{V}$ such that $f_{\mathcal{V}}(U)\subseteq V$. Since $X$ is connected, there is a chain $U_1, U_2,\ldots ,U_n$ in $\mathcal{U}$ that joins $x$ and $y$, therefore $V_1\cap f_{\mathcal{V}}(X), V_2\cap f_{\mathcal{V}}(X),\ldots ,V_n\cap f_{\mathcal{V}}(X)$, where $f_{\mathcal{V}}(U_i)\subseteq V_i\in \mathcal{V}, \; \forall i\in \{ 1,2,\ldots ,n \}$, is a chain in $\mathcal{V}$.
\end{proof}

This theorem generalized the fact that if $X$ is connected and there exists a continuous surjection $f:X\to Y$, then $Y$ must be connected.

The next example shows a case where the inheritance of connectedness cannot be applied by a continuous function, but Theorem \ref{lasttheorem} can be applied.

\begin{example}
Let the sets $X=Y=\{ a,b,c \}$ be equipped with the topologies $\tau_X=\{ \emptyset , \{ a\} , \{ b\} , \{ a,b\} , X \}$ and $\tau_Y=\{ \emptyset , \{ c\} , \{ b,c\} , \{ a,c\} , Y \}$ respectively.

Then there is no continuous surjection between them. Indeed, if there was one, since $X$ and $Y$ have the same cardinality, it would have been a bijection, therefore $X$ and $Y$ would be homeomorphic which is not true (notice that $X$ has two open singletons, while $Y$ has only one).

However, any cover of $X$ must contain $X$ itself and any cover of $Y$ must contain $Y$. Then, by Remark \ref{rem}, any function $f:X\to Y$ (or $f:Y\to X$) is $\mathcal{V}$-continuous for any covering $\mathcal{V}$ of $Y$ (or of $X$). Hence, by taking the identity mapping for any covering $\mathcal{V}$, Theorem \ref{lasttheorem} can be applied.
\end{example}

\section{Conclusion}

Coverings have been very useful tool for defining or characterizing topological notions. Beside others, coverings were used to define $\mathcal{V}$-continuity and chain connectedness, the two main notions in this article. This allowed us to use the notion of chain connectedness to prove that connectedness can be inherited only by showing the existence of continuity up to a covering functions, for every covering of the space. By an example we showed that this is an appropriate generalization of the inheritance of connectedness by using continuity.

\end{document}